\documentclass[11pt, reqno]{amsart}
\usepackage{amssymb, amsmath, amsthm, amsfonts}
\usepackage[alphabetic]{amsrefs}
\usepackage{verbatim}
\usepackage{graphics}

\newcommand{\demo}{\begin{proof}}
\newcommand{\edemo}{\end{proof}}
\newcommand{\demoname}[1]{\begin{proof}[#1]}
\newcommand{\edemoname}{\end{proof}}
\newcommand{\demonb}{\begin{proof}}
\newcommand{\edemonb}{\renewcommand{\qedsymbol}{}\end{proof}}
\newcommand{\demonamenb}[1]{\begin{proof}[#1]}
\newcommand{\edemonamenb}{\renewcommand{\qedsymbol}{}\end{proof}}

\newcommand{\stepname}{Step}

\theoremstyle{plain}
\newtheorem{theorem}{Theorem}[section]

\newtheorem{corollary}[theorem]{Corollary}
\newtheorem{lemma}[theorem]{Lemma}

\theoremstyle{definition}
\newtheorem{example}[theorem]{Example}
\newtheorem*{definition}{Definition}
\newtheorem{step}{Step}
\newtheorem{casee}{Case}
\newtheorem{stepn}{\stepname}
\newtheorem*{stepnn}{\stepname}

\newcommand{\thm}{\begin{theorem}}
\newcommand{\ethm}{\end{theorem}}

\newcommand{\expl}{\begin{example}}
\newcommand{\eexpl}{\qex\end{example}}
\newcommand{\explnb}{\begin{example}}
\newcommand{\eexplnb}{\end{example}}

\newcommand{\defn}{\begin{definition}}
\newcommand{\edefn}{\qef\end{definition}}
\newcommand{\defnnb}{\begin{definition}}
\newcommand{\edefnnb}{\end{definition}}
\newcommand{\stp}{\begin{step}}
\newcommand{\estp}{\end{step}}
\newcommand{\estpp}{\qed\end{step}}
\newcommand{\cse}{\begin{casee}}
\newcommand{\ecse}{\end{casee}}
\newcommand{\stpn}[1]{\renewcommand{\stepname}{#1}\begin{stepn}}
\newcommand{\estpn}{\end{stepn}}
\newcommand{\estpnp}{\qed\end{stepn}}
\newcommand{\stpnn}[1]{\renewcommand{\stepname}{#1}\begin{stepnn}}
\newcommand{\estpnn}{\end{stepnn}}

\newcommand{\coro}{\begin{corollary}}
\newcommand{\ecoro}{\end{corollary}}
\newcommand{\lem}{\begin{lemma}}
\newcommand{\elem}{\end{lemma}}

\providecommand{\qexsymbol}{$\lozenge$}%
\newcommand{\mathqex}{\quad\hbox{\qexsymbol}}
\DeclareRobustCommand{\qex}{%
  \ifmmode \mathqex
  \else
    \leavevmode\unskip\penalty9999 \hbox{}\nobreak\hfill
    \quad\hbox{\qexsymbol}%
  \fi
}
\providecommand{\qefsymbol}{$\triangle$}%
\newcommand{\mathqef}{\quad\hbox{\qefsymbol}}
\DeclareRobustCommand{\qef}{%
  \ifmmode \mathqef
  \else
    \leavevmode\unskip\penalty9999 \hbox{}\nobreak\hfill
    \quad\hbox{\qefsymbol}%
  \fi
}

\newcommand{\enum}{\begin{enumerate}}
\newcommand{\eenum}{\end{enumerate}}

\newcommand{\nn}{\mathbb{N}}

\newcommand{\znn}{\nn \cup \{0\}}
\newcommand{\rrr}[1]{{\mathbb{R}}^{#1}}
\newcommand{\rr}{\mathbb{R}}
\newcommand{\func}[3]{{#1} \colon {#2} \to {#3}}

\newcommand{\nd}[1]{$#1$\nobreakdash-\hspace{0pt}}

\newcommand{\nsimp}[1]{\nd{#1}simplex}

\newcommand{\ncell}[1]{\nd{#1}cell}

\newcommand{\covered}{\prec}

\newcommand{\coveredeq}{\preccurlyeq}

\newcommand{\ncovered}{\nprec}

\newcommand{\fcrit}{discrete-critical}
\newcommand{\ford}{discrete-ordinary}
\newcommand{\bcrit}{polyhedral-critical}
\newcommand{\bord}{polyhedral-ordinary}
\newcommand{\fcrito}{critical in the sense of Forman}

\newcommand{\bcrito}{critical in the sense of Banchoff}

\newcommand{\banind}[2]{a({#1}, {#2})}
\newcommand{\ordcom}[1]{\Delta({#1})}

\newcommand{\dmfs}{discrete Morse functions}

\newcommand{\uptr}{short-up-troubled}
\newcommand{\uptrl}{up-troubled}
\newcommand{\dntr}{short-down-troubled}
\newcommand{\dntrl}{down-troubled}
\newcommand{\gener}{general}

\newcommand{\twide}{\nd{2}wide}

\newcommand{\dneul}{downward Eulerian}

\newcommand{\rfu}{rank function}

\newcommand{\parfu}{parity rank function}

\newcommand{\gmod}{good variation}

\newcommand{\chn}[1]{\text{chains}({#1})}
\newcommand{\norc}[1]{l(#1)}
\newcommand{\pbe}[2]{{#1}_{<{#2}}}
\newcommand{\pbee}[2]{{#1}_{\le{#2}}}

\newcommand{\cchni}[2]{ch({#1};\, {#2})}
\newcommand{\cchnii}[3]{ch({#1};\, {#2},\, {#3})}

\newcommand{\grdd}{graded}

\newcommand{\pgrdd}{parity-graded}

\newcommand{\setcar}[1]{|{#1}|}
\newcommand{\osetcar}[1]{|}
\newcommand{\csetcar}[1]{|}

\begin{document}

%\title{Discrete Morse Functions on Regular CW Complexes are Polyhedral\\ $\ $ \\ %PRELIMINARY DRAFT}
\title{Polyhedral Representation of Discrete Morse Functions on Regular CW Complexes and Posets}
\author{Ethan D.\ Bloch}
\address{Bard College\\
Annandale-on-Hudson, NY 12504\\
U.S.A.}
\email{bloch@bard.edu}
\thanks{I would like to thank Tom Banchoff for his encouragement; Sam Hsaio for his help with posets; and the Einstein Institute of Mathematics at the Hebrew University of Jerusalem, and especially Prof. Emanuel Farjoun, for their very kind hospitality during a sabbatical when parts of this paper were written.}
\date{}
%57Q99 is none of the above in PL Topology
%52B99 is none of the above in Polytopes and Polyhedra
%57R70 is critical points and critical submanifolds in Differential Topology
%06A06 is partial order, general in Ordered Sets

\begin{abstract}
It is proved that every discrete Morse function in the sense of Forman on a finite regular CW complex can be represented by a polyhedral Morse function in the sense of Banchoff on an appropriate embedding in Euclidean space of the barycentric subdivision of the CW complex; such a representation preserves critical points.  The proof is stated in terms of discrete Morse functions on a class of posets that is slightly broader than the class of face posets of finite regular CW complexes.
\end{abstract}
\maketitle

\markright{POLYHEDRAL REPRESENTATION OF DISCRETE MORSE FUNCTIONS}

\section{Introduction and Statement of the Result}
\label{secINT}

In its classical form, Morse theory is concerned with smooth functions on smooth manifolds.  See \cite{MI2} for the basics of Morse theory.  In addition to the traditional smooth approach, there are a number of discrete analogs of Morse theory, two of which have been widely cited in the literature: the first, due to Banchoff, is found in \cite{BA1}, \cite{BA2} and \cite{BA3}; the second, due to Forman, is more recent, and is found in \cite{FORM3}, \cite{FORM4} and subsequent papers.

Although Banchoff's and Forman's approaches are both widely cited, there does not appear to be in the literature a thorough discussion of the relation between these two approaches.  Such a lack of discussion is perhaps due to the the fact that upon first encounter, the two approaches appear to be quite different.  Banchoff considers finite polyhedra embedded in Euclidean space, whereas Forman considers CW complexes (not necessarily embedded).  A ``Morse function'' for Banchoff is a projection onto a straight line in Euclidean space, whereas a ``Morse function'' for Forman, called a ``discrete Morse function,'' is a map that assigns a number to each cell of a CW complex, subject to certain conditions.  

Given a projection map, Banchoff defines an index at each vertices of a polyhedron, but at no other cells, whereas given a discrete Morse function, Forman defines an index for each critical cell, which could be of any dimension.  Banchoff does not define the concept of critical vertices vs.\ ordinary vertices in \cite{BA1}, though he does do so for polyhedral surfaces in in \cite{BA2}, and we will use that approach for all dimensions.  For Forman, the distinction between critical cells vs.\ ordinary cells is of great importance.  Finally, Banchoff focuses on relating the index at vertices to polyhedral curvature, whereas Forman focuses on using critical points for the purpose of reconstructing the CW complex up to homotopy type by attaching cells.

In spite of these apparent differences, the purpose of this note is to prove that for finite regular CW complexes there is a very concrete relation between the approaches of Banchoff and Forman, as given in the following theorem.  This theorem says that information about critical cells in the sense of Forman can be obtained by Banchoff's method for an appropriate embedding of the barycentric subdivision of the original finite regular CW complex.  We note that in Forman's method a critical \ncell p\ always has index $p$, and so the only question to be asked is whether or not a cell is critical, not what its index is.

We assume that the reader is familiar with Banchoff's approach as in \cite{BA1}, and Forman's approach as in \cite{FORM3}, though we will make use of only the first few sections of the latter paper.  We need the following clarification of Banchoff's method, which is taken from \cite{BA2}.  Let $K$ be a simplicial complex in some $\rrr{m}$, and let $\xi \in S^{m-1}$ be a unit vector.  In order to define the index at each vertex of $K$, Banchoff assumes that projection onto the line spanned by $\xi$ is general for $K$, which means that it yields distinct values for any two vertices of $K$ that are joined by an edge; this condition is true for almost all values of $\xi$.  With that assumption, Banchoff defines an index with respect to the projection, denoted $\banind v{\xi}$, for each vertex $v$ of $K$.  As mentioned above, Banchoff does not define the terms critical vertex vs.\ ordinary vertex in \cite{BA1}.  In \cite{BA2}, however, which treats only surfaces, he defines a vertex to be critical if and only if $\banind v{\xi} \ne 0$, and we will take that definition as the correct one for higher dimensions as well.

\thm\label{thmAE}
Let $X$ be a finite regular CW complex, and let $f$ be a discrete Morse function on $X$.  Then for any sufficiently large $m \in \nn$, and for any unit vector $\xi \in S^{m-1}$, there is a polyhedral embedding of the barycentric subdivision of $X$ in $\rrr{m}$ such that a cell in $X$ is \fcrito\ with respect to $f$ if and only if its barycenter is \bcrito\ with respect to projection onto the line spanned by $\xi$.
\ethm

For the sake of brevity, we will say ``\fcrit'' when we mean ``\fcrito,'' and ``\bcrit'' when we mean ``\bcrito,'' and similarly for ordinary cells and vertices.  

We take the barycentric subdivision of the CW complex in Theorem~\ref{thmAE} for the following reasons.  First, whereas Forman's method determines whether every cell is critical or ordinary, Banchoff's method assigns such information only to the vertices, and by taking the barycentric subdivision we obtain a single vertex corresponding to each original cell.  Second, the barycentric subdivision of a regular CW complex is a simplicial complex, and simplicial complexes are easier to embed in Euclidean space than more general cell complexes.  Third, even if the original CW complex were a simplicial complex, we would still need to take its barycentric subdivision \textit{prior to} embedding the complex in Euclidean space, because of the following simple example.  Let $K$ be a triangle together with its faces, which is a simplicial complex.  The function that assigns to teach face of the triangle its dimension is a discrete Morse function, and every face is \fcrit, as mentioned in \cite{FORM3}*{p.\ 108}.  However, for any embedding of the triangle in Euclidean space prior to barycentric subdivision, it is seen that the projection onto any appropriate line in the Euclidean space takes any point in the interior of an edge to a value lower than one of its vertices, and any such point, when taken as a vertex of a subdivision of the edge, is \bord.  Hence, if we want to recover the \fcrit\ cells by projection onto a line in Euclidean space, we need the flexibility of first taking the barycentric subdivision prior to embedding.

Although the motivation for this note was as stated above, it turns out that the only property of regular CW complexes that is needed for the proof of Theorem~\ref{thmAE} is the fact that the set of cells of a regular CW complex form a graded poset (partially ordered set) in a natural way, and that such a poset has various nice properties.  It is therefore more clear, and slightly more general, to formulate and prove our theorem in the context of posets.  

We assume that the reader is familiar with basic properties of posets.  See \cite{STAN}*{Chapter~3} for details.  All posets are assumed to be finite.  Let $P$ be a poset.  We let $<$ denote the partial order relation on $P$, and we write $a \covered b$ if $b$ covers $a$, and $a \coveredeq b$ if $a \covered b$ or $a = b$.  If $a \in P$, we let 
\[
\pbe Pa = \{x \in P \mid x < a\} \quad \text{and} \quad \pbee Pa = \{x \in P \mid x \le a\}.
\]

The order complex of $P$, denoted $\ordcom P$, is the simplicial complex with a vertex for each element of $P$, and a simplex for each non-empty chain of elements of $P$; it is a standard fact that such a construction yields a simplicial complex.  In particular, it is always possible to compute the Euler characteristic $\chi(\ordcom P)$.  If $C \subseteq P$ is a chain (always assumed non-empty), we let $\norc C$ denote the length of the chain, which is one less than the number of elements in the chain.

A function $\func {\rho}{P}{\{0, 1, \ldots, r\}}$ for some $r \in \znn$ is a \rfu\ for $P$ if it satisfies the following conditions: for $a, b \in P$, if $a$ is a minimal element then $\rho (a) = 0$, and if $a \covered b$ then $\rho (a) + 1 = \rho (b)$.    A poset is \grdd\ if it has a \rfu.  (There is some variation in the literature regarding the term ``\grdd''; for example, the definition used in \cite{STAN}*{Chapter~3} is more stringent.)

We will also need the following properties of posets, one of which is the mod $2$ version of \grdd\ posets, where we partly follow the terminology of \cite{BRAN}*{p.\ 6}.

\defnnb\label{defnAAQ}
Let $P$ be a finite poset.
\begin{enumerate}
\item\label{defnAAQ1}
The poset $P$ is \textbf{\twide} if for any $a, b, c \in P$ such that $a \covered b \covered c$, there is some $d \in P$ such that $d \ne b$ and $a \covered d \covered c$.  

\item\label{defnAAQ3}
Let $\func {\mu}P{\{0, 1\}}$ be a function.  The function $\mu$ is a \textbf{\parfu} if it satisfies the following conditions: for $a, b \in P$, if $a$ is a minimal element then $\mu (a) = 0$, and if $a \covered b$ then $1 - \mu (a) = \mu (b)$.  A poset is \textbf{\pgrdd} if it has a \parfu.

\item\label{defnAAQ2}
Let $\func {\mu}P{\{0, 1\}}$ be a \parfu.  The poset $P$ is \textbf{\dneul} if $a \in P$ and $a$ not minimal imply $\chi(\ordcom {\pbe Pa}) = (-1)^{\mu (a)+1} + 1$.\qef
\end{enumerate}
\edefnnb

We note that if a finite poset has a \rfu, then it is unique, and similarly for a \parfu.

Let $X$ be a regular CW complex.  Then the face poset of $X$, denoted $P(X)$, is the poset that has one element for each cell of $X$, where the order relation is given by $\sigma < \tau$ if $\sigma$ is in the boudary of $\tau$, for cells $\sigma$ and $\tau$ of $X$.  The poset $P(X)$ is ranked, where the rank of a cell in $X$ is its dimension.  It is a standard fact that $\ordcom {P(X)}$ and $X$ have homeomorphic underlying spaces.  The topological name for $\ordcom {P(X)}$ is the barycentric subdivision of $X$; if $X$ is a simplicial complex, then $\ordcom {P(X)}$ is combinatorially the same as the usual barycentric subdivision of $X$.  See \cite{L-W} or \cite{BJOR1} for details.  The poset $P(X)$ is \twide\ by \cite{FORM3}*{Theorem~1.2}.  The function that assigns each cell of $X$ the number $0$ or $1$ depending upon whether the dimension of the cell is even or odd is clearly a \parfu\ on $P(X)$.  The poset $P(X)$ is \dneul, because for each $\sigma \in P(X)$, the interval $\pbe {P(X)}{\sigma}$ is the set of all cells in the boundary of $\sigma$, which is a sphere, and hence has the appropriate Euler characteristic.

\begin{figure}[htbp]
\centering\includegraphics{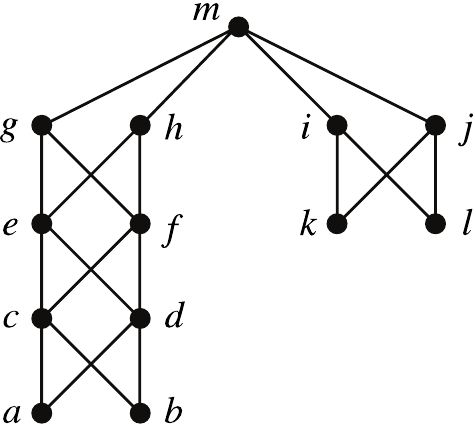}
\caption{}
\label{figMCOM2}
\end{figure}

Although the face poset of a regular CW complex is \twide, has a \parfu, and is \dneul, not every poset satisfying these three properties is the face poset of a regular CW complex.  For example, let $P$ be the poset shown in Figure~\ref{figMCOM2}.  The reader may verify that the three properties hold for $P$.  However, the poset $P$ is not the face poset of a regular CW complex, because if it were, then the interval $\pbe Pm$ would be the face poset of the boundary of cell $m$, and hence $\ordcom {\pbe Pm}$ would be a sphere, and yet $\ordcom {\pbe Pm}$ is not connected.

The original definition of discrete Morse functions in \cite{FORM3} is for discrete Morse functions on CW complexes, but the same definition can be used without change for discrete Morse functions on posets.

\defnnb\label{defnAAY}
Let $P$ be a poset, and let $\func fP{\rr}$ be a function.  
\begin{enumerate}
\item
The map $f$ is a \textbf{discrete Morse function} if the following condition holds: for each $b \in P$, there is at most one $a \in P$ such that $a \covered b$ and $f(a) \ge f(b)$, and there is at most one $c \in P$ such that $b \covered c$ and $f(b) \ge f(c)$.

\item
Suppose $f$ is a discrete Morse function.  An element $b \in P$ is \textbf{\fcrit} with respect to $f$ if there is no $a \in P$ such that $a \covered b$ and $f(a) \ge f(b)$, and there is no $c \in P$ such that $b \covered c$ and $f(b) \ge f(c)$; otherwise $b$ is \textbf{\ford} with respect to $f$.\qef
\end{enumerate}
\edefnnb

The following lemma is a restatement for posets of Lemma~2.5 of \cite{FORM3}; the original proofs works for posets, so we do not give it here.

\lem\label{lemAAZ}
Let $P$ be a finite poset, and let $\func fP{\rr}$ be a discrete Morse function.  Suppose that $P$ is \twide.  If $b \in P$, there cannot be both some $a \in P$ such that $a \covered b$ and $f(a) \ge f(b)$, and some $c \in P$ such that $b \covered c$ and $f(b) \ge f(c)$.
\elem

Lemma~\ref{lemAAZ} is not true if the assumption that $P$ is \twide\ is dropped.  For example, let $P = \{0, 1, 2\}$ have the usual total order, and let $\func fP{\rr}$ be defined by $f(x) = 2 - x$ for $x \in P$.  Then $f$ is a discrete Morse function on $P$, but it does not satisfy the conclusion of the lemma.

Our main theorem, to be proved in Section~\ref{secPRO}, is the following.

\thm\label{thmAEP}
Let $P$ be a finite poset, and let $\func fP{\rr}$ be a discrete Morse function.  Suppose that $P$ is \twide, is \pgrdd\ with \parfu\ $\mu$, and is \dneul.  Then for any sufficiently large $m \in \nn$, and for any unit vector $\xi \in S^{m-1}$, there is a polyhedral embedding $\func {\phi}{\ordcom P}{\rrr{m}}$ such that the projection of $\rrr{m}$ onto the line spanned by $\xi$ is general for $\phi(\ordcom P)$, and such that for every $b \in P$, the index of $\phi(b)$ with respect to this projection is given by
\begin{equation}\label{eqAFA}
\banind {\phi(b)}{\xi} = 
\begin{cases}
(-1)^{\mu (b)},&\text{if $b$ is \fcrit\ with respect to $f$}\\
0,&\text{if $b$ is \ford\ with respect to $f$.}
\end{cases}
\end{equation} 
\ethm

By the properties of the face poset of a CW complex stated above, and using the definition of \bcrit, Theorem~\ref{thmAE} is an immediate corollary of Theorem~\ref{thmAEP}.

Another consequence of Theorem~\ref{thmAEP} is the following.  In \cite{FORM3}*{Section~3}, analogs of some fundamental results for smooth Morse functions are proved for discrete Morse functions on regular CW complexes.  One such result is \cite{FORM3}*{Corollary 3.7 (ii)},which is part of the Weak Morse Inequalities.  We can use Theorem~\ref{thmAEP}, together with \cite{BA1}*{Theorem~1}, to deduce the following analogous result for posets.

\coro\label{coroAEQ}
Let $P$ be a finite poset, and let $\func fP{\rr}$ be a discrete Morse function.  Suppose that $P$ is \twide, is \pgrdd\ with \parfu\ $\mu$, and is \dneul.  For each $i \in \{0, 1\}$, let
\[
N_i = \setcar {\{x \in P \mid \mu(x) = i, \text{ and } x \text { is \fcrit\ with respect to } f\}}.
\]
Then $N_0 - N_1 = \chi(\ordcom P)$.
\ecoro

\demo
By Theorem~\ref{thmAEP}, there is some $m \in \nn$, some unit vector $\xi \in S^{m-1}$, and a polyhedral embedding $\func {\phi}{\ordcom P}{\rrr{m}}$ such that the projection of $\rrr{m}$ onto the line spanned by $\xi$ is general for $\phi(\ordcom P)$, and such that Equation~\ref{eqAFA} holds.

For each $i \in \{0, 1\}$, let
\[
B_i = \{x \in P \mid \mu(x) = i, \text{ and } x \text { is \fcrit\ with respect to } f\},
\]
and let
\[
D = \{x \in P \mid x \text { is \ford\ with respect to } f\}.
\]
Then
\begin{align*} 
N_0 - N_1 &= \setcar {B_0} - \setcar {B_1} = \sum_{b \in B_0} (-1)^{\mu (b)} + \sum_{b \in B_1} (-1)^{\mu (b)} + \sum_{b \in D} 0\\
 &= \sum_{b \in P} \banind {\phi(b)}{\xi} = \chi(\ordcom P),
\end{align*}
where the third equality is by Equation~\ref{eqAFA}, and the fourth equality is by \cite{BA1}*{Theorem~1}.  
\edemo

We can now use Corollary~\ref{coroAEQ} to produce a new proof of \cite{FORM3}*{Corollary 3.7 (ii)}, which is restated in the following corollary.

\coro\label{coroAER}
Let $X$ be a finite regular \nd{r}dimensional CW complex, and let $f$ be a discrete Morse function on $X$.  For each $i \in \{0, \ldots, r\}$, let $M_i$ be the number of \nd{i}cells of $X$ that are \fcrit\ with respect to $f$.  Then $\sum_{i=0}^r (-1)^iM_i = \chi(X)$.
\ecoro

\demo
As noted earlier, the face poset $P(X)$ satisfies the hypothesis of Corollary~\ref{coroAEQ}.  Because the \parfu\ on $P(X)$ is the result of assigning to each in $X$ its dimension mod two, it is seen that $\sum_{i=0}^r (-1)^iM_i = N_0 - N_1$.  As mentioned previously, it is a standard fact that $\ordcom {P(X)}$ and $X$ have homeomorphic underlying spaces.  Hence $\chi(\ordcom {P(X)}) = \chi(X)$.  The desired result now follows immediately from Corollary~\ref{coroAEQ}.
\edemo

Our proof of Corollary~\ref{coroAER} is not shorter or simpler than the proof in \cite{FORM3}*{Corollary 3.7 (ii)}, because our proof relies upon the proof of Theorem~\ref{thmAEP}, but our proof is nonetheless quite different from the original proof in that the former is essentially combinatorial, whereas the latter uses topological concepts such as homotopy equivalence.

\section{Proof of the Main Theorem}
\label{secPRO}

We start with the following three lemmas, the first two of which are very simple, and the third of which is the bulk of our work.

\lem\label{lemAB}
Let $V$ be a finite set with $n$ elements, where $n \ge 1$, and let $\func fV{\rr}$ be a function.  Then there is a map $\func {\psi}V{\rrr{n}}$ such that $\psi(V)$ spans an \nsimp {(n-1)}, and that for each vertex $v \in V$, the projection of $\psi(v)$ onto the \nd{x}axis equals $f(v)$.
\elem

\demo
The proof is by induction on $n$.  If $n = 1$, let $v$ be the single element of $V$, and then define $\psi(v) \in \rr$ to be $\psi (v) = f(v)$.  Now suppose the result is true for $n - 1$, where $n \ge 2$.  Let $w \in V$, and let $V' = V - \{w\}$.  Because $V'$ has at least one element, then by the inductive hypothesis there is a map $\func {\phi}{V'}{\rrr{n-1}}$ such that $\phi(V')$ spans an \nsimp {(n-2)}, and that for each vertex $v \in V'$, the projection of $\phi(v)$ onto the \nd{x}axis equals $f(v)$.  We can think of $\rrr{n-1}$ as sitting in $\rrr{n}$ in the usual way, and hence we can think of $\phi$ as a map $V' \to \rrr{n}$.  Let $\func {\psi}V{\rrr{n}}$ be defined by letting $\psi|_{V'} = \phi$, and letting $\psi(w)$ be a point in $\rrr{n}$ with first coordinate equal to $f(w)$, and last coordinate not equal to zero.  Because $\psi(w)$ can be joined to $\psi(V')$, we see that $\psi(V)$ spans an \nsimp {(n-1)}, and it is evident by definition that for each vertex $v \in V$, the projection of $\psi(v)$ onto the \nd{x}axis equals $f(v)$.
\edemo

For the next lemma, we need the following notation.  Let $P$ be a poset.  If $S \subseteq P$, we let $\chn S$ denote the set of non-empty chains in $S$.  If $b, s, t \in P$, and if $s \coveredeq b$ and $t \coveredeq b$, we let
\begin{align*}
%\cchni bs &= \{C \in \cchn b \mid s \in C\}\\
%\cchnii bst &= \{C \in \cchn b \mid s \in C \text{ and } t \in C\}\\
%\cchnii b{\lnot s}t &= \{C \in \cchn b \mid s \notin C \text{ and } t \in C\}.
\cchni bs &= \{C \in \chn {\pbee Pb} \mid s \in C\}\\
\cchnii bst &= \{C \in \chn {\pbee Pb} \mid s \in C \text{ and } t \in C\}\\
%\cchnii b{\lnot s}t &= \{C \in \chn {\pbee Pb} \mid s \notin C \text{ and } t \in C\}.
\cchnii b{\lnot s}t &= \{C \in \chn {\pbee Pb} \mid s \notin C \text{ and } t \in C\}.
\end{align*}

\lem\label{lemAAR}
Let $P$ be a poset.  Suppose that $P$ is \twide, is \pgrdd\ with \parfu\ $\mu$, and is \dneul.  Let $a, b \in P$.
\begin{enumerate}
\item\label{lemAAR1}
$\displaystyle \sum_{C \in \cchni bb} (-1)^{\norc C} = (-1)^{\mu (b)}$.

\item\label{lemAAR2}
If $a \covered b$, then $\displaystyle \sum_{C \in \cchnii b{\lnot a}b} (-1)^{\norc C} = 0$.

\item\label{lemAAR3}
If $a \covered b$, then $\displaystyle \sum_{C \in \cchni ba} (-1)^{\norc C} = 0$.
\end{enumerate}
\elem

\demo
For Part~(\ref{lemAAR1}), there are two cases.  First, suppose that $b$ is a minimal element of $P$.  Therefore $\mu (b) = 0$.  Also, we see that $\cchni bb = \{\{b\}\}$, and therefore $\sum_{C \in \cchni bb} (-1)^{\norc C} = (-1)^0 = (-1)^{\mu (b)}$.  Second, suppose $b$ is not a minimal element.  There is a bijective map from $\cchni bb - \{\{b\}\}$ to $\chn {\pbe Pb}$, where the map is obtained by taking each chain in the former set and removing $b$.  This map shortens the length of each chain by $1$.  Using the definition of the order complex together with the definition of \dneul, we have
\begin{align*}
\sum_{C \in \cchni bb} (-1)^{\norc C} \quad &= \sum_{D \in \chn {\pbe Pb}} (-1)^{\norc D + 1} + (-1)^{\norc {\{b\}}}\\
 &= -\sum_{D \in \chn {\pbe Pb}} (-1)^{\norc D} + (-1)^0 = - \chi(\ordcom {\pbe Pb}) + 1\\
 &= -\left[(-1)^{\mu (b)+1} + 1\right] + 1 = (-1)^{\mu (b)}.
\end{align*}

For Part~(\ref{lemAAR2}), we observe that 
\[
\cchnii b{\lnot a}b = \cchni bb - \cchnii bab.
\]
There is a bijective map from $\cchnii bab$ to $\cchni aa$, where the map is obtained by taking each chain in the former set and removing $b$.  This map shortens the length of each chain by $1$.  We then use Part~(\ref{lemAAR1}), together with the fact that $\mu (a) = 1 - \mu (b)$, to see that
\begin{align*}
\sum_{C \in \cchnii b{\lnot a}b} (-1)^{\norc C} &= \sum_{C \in \cchni bb} (-1)^{\norc C} - \sum_{D \in \cchnii bab} (-1)^{\norc D}\\
 &= \sum_{C \in \cchni bb} (-1)^{\norc C} - \sum_{D \in \cchni aa} (-1)^{\norc D + 1}\\
&= (-1)^{\mu (b)} - \left[-(-1)^{\mu (a)}\right] = 0.
\end{align*}

The proof of Part~(\ref{lemAAR3}) is similar to the proof of Part~(\ref{lemAAR2}), and we omit the details.
\edemo

\lem\label{lemAD}
Let $P$ be a poset, and let $\func fP{\rr}$ be a discrete Morse function.  Suppose that $P$ is \twide.  Then there is a discrete Morse function $g$ on $P$ that satisfies the following properties.  Let $x, y, z, w \in P$.
\enum
\item\label{lemAD1}
An element of $P$ is \fcrit\ with respect to $f$ if and only it is \fcrit\ with respect to $g$.

\item\label{lemAD3}
If $x \ne y$, then $g(x) \ne g(y)$

\item\label{lemAD2}
If $z < x \covered y < w$ and $g(x) < g(y)$, then $g(z) < g(y)$ and $g(x) < g(w)$.
\eenum
\elem

\demo
Let $a \in P$.  We say that $a$ is \textbf{\uptrl} (respectively \textbf{\uptr}) with respect to $f$ if there are $x, y \in P$ such that $a < x \covered y$ (respectively $a \covered x \covered y$) and that $f(x) < f(y) \le f(a)$.  We say that $a$ is \textbf{\dntrl} (respectively \textbf{\dntr}) with respect to $f$ if there are $z, w \in P$ such that $w \covered z < a$ (respectively $w \covered z \covered a$) and that $f(a) \le f(w) < f(z)$.  

Let $\func hP{\rr}$ be a function.  We say that $h$ is a \textbf{\gmod} of $f$ if $h$ is a discrete Morse function, and if an element of $P$ is \fcrit\ with respect to $f$ if and only if it is \fcrit\ with respect to $h$.

Step 1:  We will define a \gmod\ of $f$ that has no \uptr\ elements.  

It is a standard result that there is a total order on the set $P$ that is consistent with the original partial order $<$ on $P$.
Suppose that such a total order has been chosen.  We proceed recursively according to the total order, modifying $f$ once for each element of $P$.

Let $a \in P$ be the least element of $P$ with respect to the total order.  Then $a$ is a minimal element with respect to $<$.  If $a$ is not \uptr, we do not modify $f$ at this stage.  Now suppose that $a$ is \uptr.  Then there are $x, y \in P$ such that $a \covered x \covered y$ and $f(x) < f(y) \le f(a)$.  Observe that $a$ and $x$ are both \ford\ with respect to $f$.  By the definition of \dmfs, we know that $f(a) < f(z)$ for all $z \in P$ such that $a \covered z$ and $z  \ne x$.  By Lemma~\ref{lemAAZ} we know that if $b \in P$ and $x \covered b$, then $f(x) < f(b)$.  We then modify $f$ by decreasing the value of $f(a)$ so that it is strictly greater than $f(x)$, and strictly less than $f(b)$ for all $b \in P$ such that $x \covered b$.  The modified $f$ is a \gmod\ of $f$, and $a$ is not \uptr\ for the modified $f$.  To avoid cumbersome notation, we use $f$ to denote the modified $f$. 

Now suppose that $f$ has been modified one element of $P$ at a time so that the resulting function is a \gmod\ of $f$, and that the first $k - 1$ elements of $P$ in the total order are not \uptr.  Let $e \in P$ be the \nd{k}th element of $P$ in the total order.  If $w \in P$ and $w < e$, then $w$ is prior to $e$ in the total order, and hence $w$ is not \uptr.  If $e$ is not \uptr, we do not modify $f$ at this stage.  Now suppose that $e$ is \uptr.  Then there are $x, y \in P$ such that $e \covered x \covered y$ and $f(x) < f(y) \le f(e)$.  As before, we know that that $e$ and $x$ are both \ford, that $f(e) < f(z)$ for all $z \in P$ such that $e \covered z$ and $z  \ne x$, and that $f(x) < f(b)$ for all $b \in P$ such that $x \covered b$.  

Suppose that there is some $h \in P$ such that $h \covered e$ and $f(x) \le f(h)$.  Because $P$ is \twide, there is some $t \in P$ such that $t \ne e$ and $h \covered t \covered x$.  Because $f(x) \le f(e)$, then by the definition of \dmfs\ we know that $f(t) < f(x)$.  Because $f(x) \le f(h)$, we deduce that $h$ is \uptr, which is a contradiction.  Hence $f(d) < f(x)$ for all $d \in P$ such that $d \covered e$.  

We now modify $f$ by decreasing the value of $f(e)$ so that it is strictly greater than $f(x)$, and strictly less than $f(b)$ for all $b \in P$ such that $x \covered b$.  By the previous paragraph, we see that the modified $f$ is a \gmod\ of $f$, and now $e$ is not \uptr; the elements of $P$ that are less than $e$ in the total order remain not \uptr.  

By recursion, we can modify $f$ so that the resulting function is a \gmod\ of $f$ that has no \uptr\ elements.

Step 2: We prove that the modified $f$ has no \uptrl\ elements.  Suppose to the contrary that there is some $a \in P$ that is \uptrl.  Then there are $x, y \in P$ such that $a < x \covered y$ and $f(x) < f(y) \le f(a)$.  Because $a$ is not \uptr, then $a \ncovered x$.  By a standard fact about finite posets, there are $b_1, b_2, \ldots, b_q \in P$, with $q \ge 1$, such that $a \covered b_1 \covered b_2 \covered \cdots  \covered b_q \covered x \covered y$.  Without loss of generality, we may assume that $a$ was chosen so that $q$ is minimal for all possible \uptrl\ elements.  This minimality implies that $f(b_j) < f(y) \le f(a)$ for all $j \in \{1, \ldots, q\}$.  It follows in particular that $f(b_1) < f(a)$.  By Lemma~\ref{lemAAZ} we see that $f(b_1) < f(b_2)$, where we replace $b_2$ with $x$ if $q = 1$.  Because $a \covered b_1 \covered b_2$ and $f(b_1) < f(b_2) < f(a)$, we deduce that $a$ is \uptr, which is a contradiction.

Step 3: We will further modify $f$ so that the resulting function is a \gmod\ of $f$ that still has no \uptrl\ elements, and will now also have no \dntr\ elements.  The modification is the same as in Step~1, except that it is upside down.  The only question is whether we can perform this modification in such a way that it does not cause any elements to become \uptr; if we can make sure that no element becomes \uptr\ as a result of this modification, then by Step~2 no element will be \uptrl.

We proceed recursively, again using the total order on $P$ given in Step~1, though this time starting from the greatest element with respect to the total order, and proceeding downward.  Let $q \in P$ be the greatest element of $P$ with respect to the total order.  Modify $f$ analogously to the way we modified $f$ at the least element $a$ in Step~1, so that the modified $f$ is a \gmod\ of $f$, and $q$ is not \dntr\ for the modified $f$.  This modification of $f$, which is done by possibly increasing the value of $f(q)$, cannot cause $q$ or any element that is less than $q$ with respect to $<$ to become \uptr, and because $q$ is a maximal element with respect to $<$, there is nothing else that could become \uptr\ as a result of this modification.

Now suppose that $f$ has been modified so that the resulting function is a \gmod\ of $f$, that the last $k - 1$ elements of $P$ in the total order are not \dntr, and that there are no \uptr\ elements.  Hence by Step~2 there are no \uptrl\ elements.  Let $u \in P$ be the \nd{k}th from last in the total order.  If $v \in P$ and $u < v$, then $v$ is after $u$ in the total order, and hence by hypothesis $v$ is not \dntr.  If $u$ is not \dntr, we do not modify $f$ at this stage.  Now suppose that $u$ is \dntr.  Then there are $z, w \in P$ such that $w \covered z \covered u$ and $f(u) \le f(w) < f(z)$.  Analogously to Step~1, we know that $f(v) < f(u)$ for all $v \in P$ such that $v \covered  u$ and $v \ne z$, that $f(h) < f(z)$ for all $h \in P$ such that $h \covered z$, and that $f(z) < f(p)$ for all $p \in P$ such that $u \covered p$.

We could proceed analogously to Step~1, and modify $f$ by increasing the value of $f(u)$ so that it is strictly less than $f(c)$, and strictly greater than $f(h)$ for all $h \in P$ such that $h \covered c$, in which case the modified $f$ would be a \gmod\ of $f$, and now $u$ would not be \dntr, and the elements of $P$ that are greater than $u$ in the total order would remain not \dntr.  However, by increasing the value of $f(u)$ in this way, we might cause $u$ to become \uptr, and so we do not modify $f$ yet, but rather make the following additional observation.  

Suppose that there are $x, y \in P$ such that $u \covered x \covered y$ and $f(x) < f(y)$.  Let $r \in P$ be such that $r \covered c$.  Then $r \covered c \covered u \covered x \covered y$, and hence $r < x \covered y$.  Because $r$ is not \uptrl, then $f(r) < f(y)$.  We now modify $f$ by increasing $f(u)$ so that it is strictly less than $f(c)$, and strictly greater than $f(h)$ for all $h \in P$ such that $h \covered c$, and strictly less than $f(i)$ for all $i \in P$ such that there is some $j \in P$ such that $u \covered j \covered i$ and $f(j) < f(i)$.  We then see that the modified $f$ is a \gmod\ of $f$, and now $u$ is not \dntr\ and not \uptr, and the elements of $P$ that are greater than $u$ in the total order remain not \dntr.  The one remaining question is whether any element of $P$ other than $u$ has become \uptr\ as a result of this modification of $f$.  The only possible elements of $P$ that could be become \uptr\ as a result of increasing $f(u)$ are elements $s \in P$ for which there exist an element $t \in P$ such that $s \covered u \covered t$ or $s \covered t \covered u$; however, it is seen that in either such case, increasing the value of $f(u)$ could not make $s$ become \uptr\ if were is not already such prior to increasing $f(u)$.

By recursion, we can modify $f$ so that the resulting function is a \gmod\ of $f$ that has no \dntr\ elements and no \uptr\ elements.

Step 4: Similarly to Step~2, it can be proved that the modified $f$ has no \dntrl\ elements, as well as no \uptrl\ elements.  Hence, the modified $f$ satisfies Parts~(\ref{lemAD1}) and (\ref{lemAD2}) of the lemma.

Step 5: Let $q \in P$.  We say that $q$ is \textbf{\gener} with respect to $f$ if $f(q) \ne f(x)$ for all $x \in P - \{q\}$.

We will even further modify $f$ so that the resulting function is a \gmod\ of $f$ that still has no \uptrl\ elements and no \dntrl\ elements, and now also has all elements \gener.

As before, we proceed recursively, using the total order on $P$ given in Step~1.  Let $a \in P$ be the least element of $P$ with respect to the total order.  We then modify $f$ by increasing $f(a)$ slightly, in such a way that $a$ is \gener\ after the modification, and that nothing in $f(P)$ is between the original value of $f(a)$ and the new value; such a modification is possible because $f(P)$ is finite.  It is then seen that if $x, y \in P$ are such that $f(x) < f(y)$ prior to the modification, it must still be the case that $f(x) < f(y)$ after the modification.  It follows that no elements of $P$ can become \uptrl\ or \dntrl\ as a result of the modification, that the modified $f$ is still a discrete Morse function, and that if an element of $P$ is \fcrit\ prior to the modification, then it remains so after the modification.  Suppose that $a$ is \ford\ prior to the modification.  Because $a$ is a minimal element with respect to $<$, then it must be the case that before the modification $f(a) \ge f(b)$ for some $b \in P$ such that $a \covered b$, and hence $a$ will continue to be \ford\ after the modification.  A similar argument shows that no other element of $P$ can change from \ford\ to \fcrit\ as a result of the modification.  Hence the modified $f$ is a \gmod\ of $f$, has no \uptr\ elements and no \dntr\ elements, and $a$ is \gener.

Now suppose that $f$ has been modified so that the resulting function is a \gmod\ of $f$, that the first $k - 1$ elements of $P$ in the total order are \gener, and that there are no \uptrl\ elements and no \dntrl\ elements.  Let $e \in P$ be the \nd{k}th element of $P$ in the total order.  As before, we modify $f$ by increasing $f(e)$ slightly, in such a way that $e$ is \gener\ after the modification, and that nothing in $f(P)$ is between the original value of $f(e)$ and the new value.  Once again $e$ is \gener\ after the modification, the modified $f$ is still a discrete Morse function, no elements of $P$ can become \uptrl\ or \dntrl\ as a result of the modification, and if an element of $P$ is \fcrit\ prior to the modification, then it remains so after the modification.  Suppose that $e$ is \ford\ prior to the modification.  First, suppose that there is some $x \in P$ such that $x \covered e$ and $f(x) \ge f(e)$ prior to the modification.  Because $x < e$, then $x$ is prior to $e$ in the total order, and hence $x$ is \gener.  Therefore $f(x) > f(e)$ prior to the modification, and this inequality will still hold after the modification, and hence $e$ will remain \ford.  Second, suppose that there is some $y \in P$ such that $e \covered y$ and $f(e) \ge f(y)$ prior to the modification.  Then $f(e) > f(y)$ after the modification, and hence $e$ will remain \ford.  As before, no other element of $P$ can change from \ford\ to \fcrit\ as a result of the modification.  Hence the modified $f$ is a \gmod\ of $f$, has no \uptr\ elements and no \dntr\ elements, and first $k$ elements of $P$ in the total order are \gener.

By recursion, we can modify $f$ so that the resulting function is a \gmod\ of $f$ that has no \uptrl\ elements and no \dntrl\ elements, and has all elements \gener, which proves the lemma.
\edemo

We are now ready for the proof of our main theorem.

\demoname{Proof of Theorem~\ref{thmAEP}}
We will show that an embedding with the desired property can be found for a single choice of $\rrr{m}$ and with respect to $\xi$ being the unit vector in the direction of the positive \nd{x}axis.  It will then follow immediately that an appropriate embedding can be found in $\rrr{k}$ for $k > m$ with respect to the same $\xi$ by using using the usual embedding of $\rrr{m}$ in $\rrr{k}$.  Appropriate embeddings with respect to any unit vector $\xi \in S^{k-1}$ can be found by rotating and translating the original embedding. 

Let $g$ be the discrete Morse function on $P$ obtained by applying Lemma~\ref{lemAD} to $f$.  By Part~(\ref{lemAD1}) of the lemma, it will suffice to prove the theorem with $f$ replaced by $g$.  

Suppose that $P$ has $k$ elements.  By Lemma~\ref{lemAB} there is a map $\func {\psi}P{\rrr{k}}$ such that $\psi(P)$ spans an \nsimp {(k-1)}, and that for each vertex $v \in V$, the projection of $\psi(v)$ onto the \nd{x}axis equals $g(v)$.  Because $\ordcom P$ is a simplicial complex with $k$ vertices, it can be identified with a subcomplex of the \nsimp {(k-1)}\ spanned by $\psi(P)$.  Hence we can think of $\psi$ as inducing a polyhedral embedding $\func {\phi}{\ordcom P}{\rrr{k}}$, where $\phi(v) = \psi(v)$ for all $v \in P$, and where we think of $P$ as the set of vertices of $\ordcom P$.

By Lemma~\ref{lemAD} (\ref{lemAD3}) we know that if $a, b \in P$ and $a \ne b$, then $g(a) \ne g(b)$.  It follows that if $\psi(a)$ and $\psi(b)$ are vertices of $\ordcom P$ that are joined by an edge, then $g(a) \ne g(b)$, and hence the projection of $\psi(a)$ onto the \nd{x}axis does not equal the projection of $\psi(b)$ onto the \nd{x}axis.  We can therefore define Banchoff's index at the vertices of $\phi(\ordcom P)$, where the projection is onto the \nd{x}axis, and hence we can apply the notion of \bcrit\ and \bord\ to these vertices.

Let $b \in P$, so that $b$ is a vertex of $\ordcom P$.  Following \cite{BA1}, we compute the index $\banind {\phi(b)}{\xi}$ as follows.  Let $T$ denote the set of all simplices of $\ordcom P$ that contain $\phi(b)$ as a vertex and for which projection onto the \nd{x}axis has maximal value at $\phi(b)$.  Then $\banind {\phi(b)}{\xi} = \sum_{s \in T} (-1)^{\dim s}$.  We can view this last formula from a different perspective.  By the definition of $\ordcom P$, every simplex of $\ordcom P$ is a non-empty chain in $P$.  The choice of $\phi$ states that the projection of $\phi(a)$ onto the \nd{x}axis equals $g(a)$ for all $a \in P$.  Hence, we see we can think of $T$ as the set of all chains in $P$ that contain $b$, and on which $g$ is maximal at $b$.  If $C$ is a chain in $P$, then the dimension of this chain when thought of as a simplex of $\ordcom P$ is equal to $\norc C$.  Therefore $\banind {\phi(b)}{\xi} = \sum_{C \in T} (-1)^{\norc C}$.  

Suppose that $b$ is \fcrit\ with respect to $g$.  Let $v \in P$ be such that $v < b$.  If $v \covered b$, then $g(v) < g(b)$ because $b$ is \fcrit\ with respect to $g$.  Now suppose $v \ncovered b$.  By a standard fact about finite posets, there is some $z \in P$ such that $v < z \covered b$.  Because $b$ is \fcrit\ with respect to $g$, then $g(z) < g(b)$.  By Lemma~\ref{lemAD} (\ref{lemAD2}) we deduce that $g(v) < g(b)$.  A similar argument shows that $g(b) < g(u)$ for any $u \in P$ such that $b < u$.  Hence, the set $T$ consists precisely of all chains in $\pbee Pb$ that contain $b$; this set is denoted $\cchni bb$.  Lemma~\ref{lemAAR} (\ref{lemAAR1}) implies that
\[
\banind {\phi(b)}{\xi} = \sum_{C \in T} (-1)^{\norc C} = \sum_{C \in \cchni bb} (-1)^{\norc C} = (-1)^{\mu (b)}.
\]

Next, suppose that $b$ is \ford\ with respect to $g$.  Then by Lemma~\ref{lemAAZ} either there is a single $h \in P$ such that $h \covered b$ and $g(h) \ge g(b)$, or there is a single $u \in P$ such that $b \covered u$ and $g(b) \ge f(u)$, but not both. 

First, suppose that there is some $h \in P$ such that $h \covered b$ and $g(h) \ge g(b)$.  Then $g(b) < g(z)$ for all $z \in P$ such that $b \covered z$.  By the same argument used above, we know that $g(b) < g(u)$ for any $u \in P$ such that $b < u$, and hence that $T \subseteq \cchni bb$. 

Let $c \in P$ be such that $c < b$ and $c \ne h$.  If $c \covered b$, then the definition of discrete Morse functions implies that $g(c) < g(b)$.  Now suppose $c \ncovered b$.  Combining basic properties of posets with the fact that $P$ is \twide, there is some $t \in P$ such that $t \ne h$ and $c < t \covered b$.  Then by the definition of discrete Morse functions we know that $g(t) < g(b)$, and using Lemma~\ref{lemAD} (\ref{lemAD2}) it follows that $g(c) < g(b)$.  

Putting the above considerations together, we see that $h$ is the only element of $\pbee Pb$ such that $g(h) > g(b)$.  Hence $T = \cchnii b{\lnot h}b$, and Lemma~\ref{lemAAR} (\ref{lemAAR2}) implies that 
\[
\banind {\phi(b)}{\xi} = \sum_{C \in \cchnii b{\lnot h}b} (-1)^{\norc C} = 0.
\]

Second, suppose that there is some $u \in P$ such that $b \covered u$ and $g(b) \ge g(u)$.  An argument similar to the previous case shows that $T = \cchni ub$.  Lemma~\ref{lemAAR} (\ref{lemAAR3}) implies that $\banind {\phi(b)}{\xi} = 0$.
\edemoname

\end{document}